\documentclass[12pt,a4paper]{article}

\usepackage[latin1]{inputenc}
\usepackage[english]{babel}
\usepackage{amsmath,amsthm,amssymb}
\usepackage{makeidx}
\usepackage{amsfonts} 
\usepackage{textcomp}
\usepackage{stmaryrd}
\usepackage[all]{xy}
\numberwithin{equation}{section}
\swapnumbers

{\theoremstyle{definition}\newtheorem{definition}{Definition}[section]

\newtheorem{notation}[definition]{Notation}
\newtheorem{remarque}[definition]{Remark}
\newtheorem{remht}[definition]{Remark on hyperlinear groups and traces}

\newtheorem{remarques}[definition]{Remarks}

}
\newtheorem{proposition}[definition]{Proposition}
\newtheorem{lemme}[definition]{Lemma}

\newtheorem{theoreme}[definition]{Theorem}
\newtheorem{corollaire}[definition]{Corollary}

\newcommand{\cf}{\textit{cf. }}
\newcommand{\ie}{\textit{i.e. }}

\newcommand{\R}{\mathbb{R}}
\newcommand{\Q}{\mathbb{Q}}
\newcommand{\N}{\mathbb{N}}
\newcommand{\Z}{\mathbb{Z}}
\newcommand{\C}{\mathbb{C}}
\newcommand{\F}{\mathbb{F}}
\newcommand{\Tr}{\mathrm{Tr}}
\newcommand{\gS}{\mathfrak{S}}
\newcommand{\cT}{\mathcal{T}}
\newcommand{\card}{\mathrm{card}}
\newcommand{\Sp}{\mathrm{Sp}}

\begin{document}
\everymath={\displaystyle}
\renewcommand{\labelitemi}{$\bullet $}
\renewcommand{\labelitemii}{$\ast $}

\begin{center}
{\Large\bf Traces on group $C^\ast$-algebras, sofic groups \\and Lück's conjecture}
\bigskip

{\sc by G\"ul Balci and Georges Skandalis}
\end{center}

{\footnotesize
\vskip 2pt Universit\'e Paris Diderot, Sorbonne Paris Cit\'e
\vskip-4pt  Sorbonne Universit\'es, UPMC Paris 06, CNRS, IMJ-PRG
\vskip-4pt  UFR de Math\'ematiques, {\sc CP} {\bf 7012} - B\^atiment Sophie Germain 
\vskip-4pt  5 rue Thomas Mann, 75205 Paris CEDEX 13, France
\vskip-5pt e-mail: balcigul@hotmail.com\ 
\vskip-5pt \ \hskip 1.1cm georges.skandalis@imj-prg.fr
}
\bigskip

\centerline{\bf Abstract}

We give an alternate proof of a Theorem of Elek and Szabo establishing L\"uck's determinant conjecture for sofic groups. Our proof is based on traces on group C*-algebras. We briefly discuss the relation with Atiyah's problem on the integrality of $L^2$-Betti numbers.

\section*{Introduction}

In \cite{FK}, B. Fuglede and R.V. Kadison introduce a determinant function $\Delta:M\to \R_+$ where $M$ is a II$_1$ factor, by setting $\Delta (x)=\exp(\tau (\ln |x|))$  (where $\tau$ is the trace of $M$). This function satisfies many of the usual properties of a determinant. 

In connexion with many problems and conjectures concerning discrete groups, W. Lück (see \cite{Lueck2}) introduced a modified determinant by setting $\Delta _+(x)=\exp(\tau (\ln_+ |x|))$ where $\ln_+(t)=0$ if $t=0$ and $\ln_+(t)=\ln t$ for $t>0$. He conjectured that for any group $G$ and any $x\in M_n(\Z G)$ we have $\Delta_+(|x|)\ge 1$.

Lück's determinant conjecture is related with many interesting problems (\cf \cite{Lueck3, Lueck2}). In particular:\begin{itemize}
\item Let us recall a problem stated first by Atiyah in the torsion free case, and extended by various authors  (\cf \cite{Lueck2, Schick1bis, Schick2} for more details): investigate the possible values of von Neumann dimensions of the kernels of elements in $M_n(\Z G)$. Validity of Lück's conjecture ensures a kind of stability of the von Neumann dimension of the kernel of an element in $M_n(\Z G)$ and allows its computation in some cases.
\item For CW complexes whose fundamental group satisfies this conjecture, one can define $L^2$ torsion (\cite{Lueck2}). In  \cite{LSW}, Lück, Sauer and Wegner prove that the $L^2$-torsion is invariant under uniform measure equivalence.
\end{itemize}

In \cite{Schick1, Schick1bis}, Schick shows that amenable groups satisfy Lück's conjecture. He shows that the class of groups satisfying Lück's conjecture is closed under taking subgroups, direct limits and inverse limits.

In \cite{ElekSzabo2}, Elek and Szabo generalize Schick's results by proving that sofic groups satisfy Lück's conjecture.
Sofic groups were introduced by Gromov \cite{Gromov} as a generalization of both amenable and profinite groups. In a sense, sofic groups are the groups that can be well approximated by finite groups, \ie that can be almost embedded into permutation groups.

\bigskip 
In this paper we give a reformulation of Lück's conjecture in terms of traces on  $C^*(\F_\infty)$ and use it to reformulate the proof of \cite{ElekSzabo2} in a somewhat more conceptual way. We hope that it may help understanding the proofs of  \cite{Schick1} and  \cite{ElekSzabo2}.

We say that a trace $\tau$ on $C^*(\F_\infty)$ satisfies Lück's condition if for all $f\in M_n(\Z \F_\infty)$, we have $\tau(\ln_+(\vert f\vert))\ge 0$.
In this sense, a group $G$ satisfies Lück's conjecture if and only if the trace $\tau_G\circ \pi$ satisfies Lück's condition, where $\tau_G $ is the canonical trace on $G$ and $\pi$ a surjective morphism $\F_\infty\to G$ (\ie a generating system of $G$). It is then easily seen (a detailed account is given in the sequel) that:\begin{description}
\item[Fact 1.]  (proposition \ref{canon}) A permutational trace \ie a trace of the form $tr\circ \sigma$ satisfies Lück's condition where $\sigma $ is a finite dimensional representation of $\F_\infty$ by permutation matrices and $tr$ is the normalized trace on matrices.
\item[Fact 2.] (proposition \ref{closed}) The set of traces satisfying Lück's condition is closed (for the weak topology).
\item[Fact 3.] (proposition \ref{propequi}) A group $G$ is sofic if and only if the associated trace $\tau_G\circ \pi$ as above is in the closure of permutational traces. 
\end{description}

Moreover, we notice that the stability condition of the von Neumann dimension established in \cite{Schick1, Schick2} is a consequence of the following fact:

\begin{description}
\item[Fact 4.] (proposition \ref{propAt}) For $a\in M_n(\Z \F_\infty)$, the map $\tau \mapsto \dim_\tau \ker a$ is continuous on the set of traces satisfying Lück's condition. 
\end{description}

Finally, we extend this result to $a\in M_n(\overline\Q \F_\infty)$. We moreover prove that, for any trace of $C^*(\F_\infty)$  in the closure of permutational traces, the dimension $\dim_\tau (\ker a)$ does not depend on the embedding $\overline\Q\subset \C$.  We deduce the following  formulation of a result of \cite{DLMSY} to sofic groups which is proved by A. Thom in (the proof of) \cite[Theorem 4.3]{Thom}.

\begin{description}
\item[Fact 5] (corollary \ref{corindep}) Let $\Gamma$ be a sofic group and $a\in M_n(\overline\Q \Gamma)$. The von Neumann dimension (with respect to the group trace of $\Gamma$)  does not depend on the embedding $\overline\Q\subset \C$.
\end{description}

\bigskip This paper is organized as follows : 
In the first section, we fix notation and recall definitions of the determinant of Fuglede-Kadison and the modified determinant of Lück.\\
In the second section, we define Lück's condition for a trace and establish facts 1 and 2 above. \\
In the third section we recall the definition of a sofic group and establish fact 3.\\
In section 4, we explain the relation with Atiyah's problem and establish facts 4 and 5.

\emph{All traces that we consider throughout the paper are positive finite traces.}

\section{Positive traces and determinants}

\subsection{Traces and semi-continuous functions}

Let $A$ be a unital $C^*$-algebra. We endow the set $\cT_A$ of (finite positive) traces on $A$ with the pointwise convergence.

Let $\tau \in \cT_A$ be a  trace on $A$ and $a$ a self-adjoint element of $A$, and $\mu_{\tau,a}$ the corresponding spectral measure. If $f:\Sp\, a\to \R\cup \{+ \infty\}$ is a lower semi-continuous function, we may write $f=\sup f_n$ where $f_n$ is an increasing sequence of continuous functions. Since $f$ is bounded below, we may define $\mu_{\tau,a}(f)\in \R\cup\{+\infty\}$ and we have, by the monotone convergence theorem,  $\mu_{\tau ,a}(f)=\sup \mu_{\tau,a}(f_n)=\sup \tau (f_n(a))$. In the sequel, this ``number" will be denoted by $\tau(f(a))$. 

In the same way, we define $\tau(f(a))\in\R\cup \{- \infty\}$ for $f:\Sp\, a\to \R\cup \{- \infty\}$ upper semi continuous.

We obviously have:

\begin{proposition}\label{propsc}
If $f:\Sp\,a\to \R\cup \{+ \infty\}$ is lower (\emph{resp.} upper) semi-continuous, then the map $\tau\mapsto \tau(f(a))$ is lower (\emph{resp.} upper) semi-continuous.
\begin{proof}
Assume $f=\sup f_n$ is lower semi-continuous where $(f_n)$ is an increasing sequence of continuous functions.
Then $\tau \mapsto \tau(f(a))$ is the supremum of the continuous functions $\tau \mapsto \tau(f_n(a))$ and is therefore lower semi-continuous.

The upper semi-continuous case is obtained by replacing $f$ by $-f$.
\end{proof}
\end{proposition}

\begin{remarque} \label{rem}
Let $\varphi : A\to B$ be a unital morphism of $C^\ast$-algebras and $\tau$ a trace on $B$.
Then, for every self-adjoint element $a\in A$ and  every lower semi-continuous function $f:\Sp\, a\to \R\cup \{+\infty\}$, writing $f=\sup f_n$ with continuous $f_n$, we find  $f_n(\varphi(a))=\varphi(f_n( a))$; hence passing to the supremum, one gets $\tau(f (\varphi(a)) = \tau\circ\varphi(f(a))$.

The same equality holds of course for upper semi-continuous $f$.
\end{remarque}

\subsection{The Fluglede-Kadison determinant (\cite{FK})}

The function $\ln:\R_+\to \R\cup \{- \infty\}$ is upper semi-continuous on $\R_+$. 

Let $A$ be a unital $C^*$-algebra and $\tau \in \cT_A$ a (finite, positive) trace on $A$.
 
The Fuglede-Kadison determinant of $x\in A$ is $\Delta_\tau(x)= \exp(\tau(\ln(|x|))) $.

Recall from \cite{FK} that, for $x,y\in A$, we have $\Delta_\tau(xy)=\Delta_\tau(x)\Delta_\tau(y)$.

\subsection{Lück's modified determinant (\cite{Lueck2})}

For $t\in\R_+$, set $$\ln_+(t)=\begin{cases} 0 &\hbox{ if } t=0\\
\ln(t) &\hbox{ if } t\neq 0 \end{cases}$$
 
 Note that $\ln_+$ is also upper semi-continuous.
 
Lück's modified  determinant of $x\in A$ is $\Delta^+_\tau(x)= \exp(\tau(\ln_+(|x|))) $.

If $\Tr_n$ is the unnormalized trace of $M_n(\C)$, then $\Delta_{\Tr_n}(a)=|\mathrm{Det}(a)|=\mathrm{Det}(|a|)$ where $\mathrm{Det}$ is the usual determinant, and $\Delta^+_{\Tr_n}(a)$ is the product of nonzero eigenvalues of $|a|$ (counted with their multiplicity).

\section{Traces satisfying Lück's condition}

If $G$ is a (discrete) group, we denote by $C^\ast(G)$ the full group $C^\ast$-algebra of $G$.

A trace on $C^*(G)$ is determined by its values on the dense subalgebra $\C G$, whence by linearity, by its value on the group elements. We will identify the set of traces on $C^*(G)$ with the set $\cT_G$ of maps $G\to \C$ which are of positive type and constant on conjugacy classes.

The  weak topology on the set $\cT_{C^*(G)}=\cT_G$ of traces on $C^*(G)$ coincides with the topology of pointwise convergence on $G$.

We denote by $\tau_G$ the canonical trace on $G$, \ie the one given by $\tau_G(x)=1$ if $x=1_G$ - the unit of $G$ and $\tau_G(x)=0$ for $x\in G\setminus\{1_G\}$.

\begin{definition}
Let $G$ be a countable group.
We say that a trace $\tau$ on $C^\ast(G)$ satisfies Lück's condition if for $n\in\N$ and $a\in M_n(\Z G)\subset C^\ast(G)\otimes M_n(\C) $, we have $(\tau\otimes \Tr_n)(\ln_+ \vert a\vert)\geq 0$ where $\Tr_n$ is the unnormalized trace on $M_n(\C)$. We denote by $\Lambda_G\subset \cT_G$ the set of traces on $C^\ast(G)$ satisfying Lück's condition\\ The group $G$ is said to satisfy Lück's conjecture if the canonical trace $\tau_G$ of $C^\ast(G)$ satisfies Lück's condition.
\end{definition}

\begin{proposition}\label{double}
Let $G$ and $H$ be groups and $\psi : G\to H$ a group homomorphism. 
We still denote by $\psi : C^\ast(G)\to C^\ast(H)$ its extension to group $C^\ast$-algebras. Let $\tau$ be a trace on $C^\ast (H)$.
\begin{enumerate}
\item If $\tau\in \Lambda_H$ (\ie $\tau $ satisfies Lück's condition),  $\tau\circ\psi\in \Lambda_G$. 
\item If $\psi$ is onto, the converse is true, i.e. if $\tau\circ\psi\in \Lambda_G$, then $\tau\in \Lambda_H$. 
\end{enumerate}

\begin{proof}
\begin{enumerate}
\item Suppose that $\tau$ satisfies Lück's condition.\\
We denote by $\psi : M_n(C^\ast (G))\to M_n(C^\ast (H))$ the extension of $\psi$ to matrices.\\
For all $a\in M_n(\Z G)$, since $\psi(a)\in M_n(\Z H)$, we have, by remark \ref{rem}, $$(\tau\otimes \Tr_n)(\psi(\ln_+ \vert a\vert))=(\tau\otimes \Tr_n)(\ln_+ \vert\psi(a)\vert)\geq 0.$$

\item  If $\psi$ is surjective, for all $a\in M_n(\Z H)$ there exists $b\in M_n(\Z G)$ such that $a=\psi(b)$.
If $\tau\circ\psi$ satisfies Lück's condition then $(\tau\otimes \Tr_n)(\ln_+(\vert a\vert))=(\tau\otimes \Tr_n)\circ\psi(\ln_+(\vert b\vert))\geq0$.
\qedhere
\end{enumerate}
\end{proof}
\end{proposition}

Let $\alpha_k :\gS_k\to\mathcal{U}_k$ be the representation of the symmetric group $\gS_k$ by permutation matrices. 
Denote also by $\alpha_k$ the associated morphism $C^\ast (\gS_k)\to M_k(\C)$. 
We denote by $\tau_k$ the trace on $C^\ast(\gS_k)$ defined by 
$\tau_k(f)=\frac{1}{k}\Tr_k(\alpha_k(f))$
for all $f\in C^\ast (\gS_k)$,  where $\Tr_k$ is the unnormalized trace on $M_k(\C)$.\\
Then, for $\sigma\in\gS_k$, we have $$\tau_k(\sigma)=\dfrac{\card\left\lbrace x\vert \sigma(x)=x\right\rbrace }{k}.$$

\begin{proposition}\label{canon} (\cf \cite{Lueck2})
Let $n\in\N^\ast$. The trace $\tau_k$ on $\gS_k$, satisfies Lück's condition.

\begin{proof}
Let $\Tr_{kn}$ be the unnormalized trace on $M_{kn}(\Z)$.
For $a\in M_n(\Z\gS_k)$ we have
$$ \tau_k\otimes \Tr_n(a)=\frac{1}{k}\Tr_{kn}(\alpha_k(a)).$$
Now, $\Delta_{\Tr_{kn}}^+(\alpha_k(|a|))^2=\exp(2k  (\tau_k\otimes \Tr_n)(\ln_+(|a|)))$ is the product of non-zero eigenvalues of $\alpha_k(a^*a)$ (counted with multiplicity): it is the modulus of the non-zero coefficient of lowest  degree of the characteristic polynomial of $\alpha_k(a^*a)$.

Since $\alpha_k(a^*a)\in M_{kn}(\Z)$, it follows that its characteristic polynomial has integer coefficients so that $\Delta^{\Tr_{kn}}_+(a^*a)\in \N^\ast$. 
\end{proof}
\end{proposition}

\begin{definition}
We call permutational trace on a group $G$, a trace of the form $\tau_k\circ f$ where $f$ is a group morphism from $G$ to $\gS_k$. 
\end{definition}

It follows from \ref{double} and \ref{canon} that permutational traces satisfy Lück's condition.

We will also use the following more general sets of traces. 
\begin{notation}
 Let $A$ be a unital $C^*$-algebra, $d\in \N$ and $a\in M_d(A)$. Let $s\in \R$. Denote by $\Lambda_{a,s}\subset \cT_A$ the set of traces on $A$ such that $(\tau\otimes \Tr_d)(\ln_+(|a|))\ge s$.
\end{notation}

\begin{proposition}\label{closed}
\begin{enumerate}
\item Let $A$ be a unital $C^*$-algebra, $d\in \N$ and $a\in M_d(A)$. Let $s\in \R$. The set $\Lambda_{a,s}$ is closed in $\cT_A$ (for the pointwise topology on traces).
\item Let $G$ be a group. The set $\Lambda_G$ is closed in the set $\cT_G$ of traces on $C^\ast(G)$ for the pointwise topology.
\end{enumerate}
\begin{proof}
\begin{enumerate}
\item Since $\ln_+$ is upper semi-continuous, the map $\tau\mapsto (\tau\otimes \Tr_d)(\ln_+|a|)$ is upper semi-continuous by prop. \ref{propsc}, therefore the set $\Lambda_{a,s}$ is closed.

\item The set $\Lambda_G$ is the intersection over all $k\in \N$ and $a\in M_d(\Z G)$ of the closed subsets $\Lambda_{a,0}$. It is closed. \qedhere
\end{enumerate}
\end{proof}
\end{proposition}

Note that the set $\Lambda_{a,s}$ only depends on the abelian $C^*$-subalgebra of $M_d(A)$ generated by $a^*a$ and that the prop. \ref{closed} can immediately be extended to all positive forms.

\section{Sofic groups and traces}

Denote by $\delta_k$ the Hamming distance on $\gS_k$ defined by
$$\delta_k(\sigma_1,\sigma_2)=\frac{1}{k}\card\left\lbrace x\in \left\lbrace 1,...,k\right\rbrace \vert \sigma_1(x)\neq\sigma_2(x)\right\rbrace $$
with $\sigma_1, \sigma_2\in\gS_k$.

\begin{remarques} \label{remarq}
\begin{enumerate}
\item The distance $\delta_k$ is left and right-invariant, \ie for $\alpha,\beta,\sigma_1,\sigma_2\in\gS_{k_n}$, we have $$\delta_k(\alpha\circ\sigma_1\circ\beta,\alpha\circ\sigma_2\circ\beta)=\delta_k(\sigma_1,\sigma_2).$$

\item For $\sigma\in\gS_{k}$, we have $\delta_k(\sigma,Id_k)=1-\tau_k(\sigma)$.
In particular, $\tau_k$ is 1-lipschitz for $\delta_k$.
\end{enumerate}
\end{remarques}

Let $(k_n)_{n\in \N}$ be a sequence of integers. Put $$\bigoplus _{n\in\N}^{c_0}\gS_{k_n}=\lbrace(\sigma_n)_{n\in\N}\in\prod_{n\in\N}\gS_{k_n} \vert \delta_{k_n}(\sigma_n,Id_{k_n})\to 0\rbrace \subset \prod_{n\in\N}\gS_{k_n}.$$
By the invariance property, it is a normal subgroup of $\prod_{n\in\N}\gS_{k_n}$.

Recall Gromov's definition of a sofic group (\cf \cite{Gromov}; see also \cite{Elekszabo1, Pestov} for a very nice presentation of sofic groups).

\begin{definition}\label{sofic}
A countable group $\Gamma$ is said to be \emph{sofic} if there exists a sequence of maps $(f_n)_{n\in\N} : \Gamma\to \gS_{k_n}$ such that : \begin{enumerate}\renewcommand \theenumi{\alph{enumi}}
\renewcommand \labelenumi{\rm {\theenumi})}
\item For all $x,y\in\Gamma$, $\delta_{k_n}(f_n(x)f_n(y), f_n(xy))\to 0$,
\item For all $x\in\Gamma$, $x\neq 1$, $\delta_{k_n}(f_n(x),Id)\to 1$.
\end{enumerate}
Such a sequence $(f_n)_{n\in\N}$ is called a \emph{sofic approximation} of $\Gamma$.
\end{definition}

\begin{remarque} Let $q : \prod_{n\in\N}\gS_{k_n} \to \prod_{n\in\N}\gS_{k_n}/\bigoplus _{n\in\N}^{c_0}\gS_{k_n}$ be the quotient map.
Property $(a)$ of definition \ref{sofic}, means that $g=q\circ (f_n)_{n\in\N} : \Gamma\to \prod_{n\in\N}\gS_{k_n}/\bigoplus _{n\in\N}^{c_0}\gS_{k_n}$ is a morphism.

In particular, if (a) is satisfied, $(f_n(1))\in \bigoplus _{n\in\N}^{c_0}\gS_{k_n}$.
Therefore conditions (a) and (b) are equivalent to (a) and\begin{description}
\item[{\rm{b')}}]  For all $x\in\Gamma$, $\tau_{k_n}(f_n(x))\to\tau_\Gamma(x)$.
\end{description}
\end{remarque}

We now prove the rather easy characterization of sofic groups in terms of traces:

\begin{proposition}\label{propequi}
Let $\Gamma$ be a countable group. Denote by $\tau_\Gamma$ the canonical trace on $C^\ast(\Gamma)$.
The following are equivalent:\begin{enumerate}\renewcommand \theenumi{\roman{enumi}}
\renewcommand \labelenumi{\rm ({\theenumi})}
\item  The group $\Gamma$ is sofic.
\item There exists a group $G$ and a surjective morphism $\varphi : G\to \Gamma$ such that $\tau_\Gamma\circ \varphi$ is in the closure of the permutational traces.
\item For every surjective morphism $\varphi : \F_\infty\to\Gamma$ (\ie for every generating system of $\Gamma$), $\tau_\Gamma\circ\varphi$ is in the closure of permutational traces.
\end{enumerate}

\begin{proof}
$(iii)\Rightarrow (ii)$ is obvious.

$(ii) \Rightarrow (i)$ Assume that there is an onto morphism $\varphi : G\to\Gamma$ and a sequence $(\pi_n)_{n\in\N}$ of morphisms $\pi_n : G\to\gS_{k_n}$ such that for $y\in G$, 
\begin{eqnarray}
\tau_{k_n}(\pi_n(y))&\to&\tau_\Gamma(\varphi(y)) \label{eqn1}
\end{eqnarray}
Let $s$ be any section of $\varphi$ and put $f_n=\pi_n\circ s$. 
We show that the sequence $f=(f_n)$ is a sofic approximation, by establishing conditions (a) and (b')

\bigskip

\begin{description}
\item[(b')] Let $x\in\Gamma$. Applying (\ref{eqn1}) to $y=s(x)$, we get $(b')$.
 
\item[(a)] The family $\pi=(\pi_n)$ determines a morphism $\pi:G\to \prod_{n\in\N}\gS_{k_n}$. For $y\in \ker \varphi$, since, $\tau_\Gamma(\varphi(y))=1$, we find $\delta_{k_n}(1,\pi_n(y))=1-\tau_{k_n}(\pi_n(y))\to 0$; therefore $\pi(y)\in \bigoplus _{n\in\N}^{c_0}\gS_{k_n}=\ker q$. It follows that $q\circ f$ is a morphism, and (a) is satisfied.
\end{description}

$(i)\Rightarrow (iii)$ Let $f=(f_n)$ be a sofic approximation of $\Gamma$. Then $q\circ f:\Gamma\to \prod_{n\in\N}\gS_{k_n}/\bigoplus _{n\in\N}^{c_0}\gS_{k_n}$ is a morphism.

Since $\F_S$ is free, the morphism $q\circ f\circ \varphi$ lifts to a morphism $\pi=(\pi_n):\F_S\to \prod_{n\in\N}\gS_{k_n}$.

Now, for $y\in \F_S$, since $\pi(y)^{-1}f\circ \varphi(y)\in \ker q$, we find $\delta_{k_n}(\pi_n(y),f_n(\varphi(y)))\to 0$, whence (since $\tau_{k_n}$ is 1-lipschitz), $|\tau_{k_n}(\pi_n(y))-\tau_{k_n}(f_n(\varphi(y)))|\to 0$. The property (b') of the sofic approximation $f$ yields $\tau_{k_n}(\pi_n(y))\to\tau_\Gamma\circ\varphi(y)$
\end{proof}
\end{proposition}

\begin{theoreme}\label{thetheo} (cf. \cite{ElekSzabo2})
Every sofic group satisfies Lück's conjecture.
\begin{proof}
Let $\Gamma$ be a sofic group.
Then, there exists a surjective morphism $\varphi : \F_\infty\to \Gamma$. By prop. \ref{propequi}, the trace $\tau_\Gamma\circ \varphi$ is the limit of permutational traces on $C^\ast(\F_\infty)$.\\
By proposition \ref{canon}, the canonical trace on $\gS_n$ satisfies Lück's conjecture and the same is true for a permutational trace on $C^\ast(\F_\infty)$ by proposition \ref{double} $(i)$.\\
Then $\tau_\Gamma\circ\varphi\in \Lambda_{\F_\infty}$ because $\Lambda_{\F_\infty}$ is closed in the set of traces on $C^\ast(\F_\infty)$ by proposition \ref{closed}.\\
Since $\varphi$ is surjective, we conclude by proposition \ref{double} $(ii)$, that $\tau_\Gamma$ satisfies Lück's conjecture. 
\end{proof}
\end{theoreme}

\begin{remarque}
A convenient way to state this result is to define the set of \emph{sofic traces} on a group $G$ as being the closure of permutational traces. Then \begin{enumerate}
\item A group is sofic if and only if the trace it defines on $\F_\infty$ is sofic. 
\item Every sofic trace satisfies Lück's condition.
\end{enumerate}

On the other hand, not every trace on $\F_\infty$ satisfies L\"uck's condition. It becomes then a quite natural question to study the set of traces on $\F_\infty$ satisfying Lück's condition, the set of sofic traces, \emph{etc.} Such a study is undertaken in \cite{Gul}.
\end{remarque}

\begin{remht}  (see \textit{e.g.} \cite{Pestov} for a definition of hyperlinear groups)
In the same way, we may define linear traces as the characters of finite dimensional representations and the set of \emph{hyperlinear traces} as the closure of the set of linear traces. One then easily shows that a group is hyperlinear if and only if the trace it defines on $\F_\infty$ is hyperlinear.

One actually proves (see  \cite{Gul} for details):

\textit{Let $\Gamma$ be a countable group. Denote by $\tau_\Gamma$ the canonical trace on $C^\ast(\Gamma)$.
The following are equivalent:\begin{enumerate}\renewcommand \theenumi{\roman{enumi}}
\renewcommand \labelenumi{\rm ({\theenumi})}
\item  The group $\Gamma$ is hyperlinear.
\item There exists a group $G$ and a surjective morphism $\varphi : G\to \Gamma$ such that the trace $\tau_\Gamma\circ \varphi$ is hyperlinear.
\item There exists a group $G$, a surjective morphism $\varphi : G\to \Gamma$ and a hyperlinear trace $\tau$ on $G$ such that $\{g\in G;\ \tau(g)=1\}=\ker \varphi$.
\item For every surjective morphism $\varphi : \F_\infty\to\Gamma$ (\ie for every generating system of $\Gamma$), the trace $\tau_\Gamma\circ\varphi$ on $\F_\infty$ is hyperlinear.
\end{enumerate}}
\end{remht}

\section{Relation with Atiyah's problem}

L\"uck's conjecture implies a kind stability of von Neumann dimension. This stability allows computing $L^2$-Betti numbers, and proving integrality in some cases (\cite{Lueck2, Lueck3, Schick2}, ...), or actually disproving their rationality (\cite{GLSZ, Austin})

\subsection{The method}

Let $A$ be a unital $C^*$-algebra , $k\in \N$ and $a\in A$. The characteristic function $\chi_0$ of $\{0\}$ is upper semi-continuous.   Given a trace $\tau$ on $A$, we put $\dim_{\tau}(\ker a)=\tau(\chi_0(|a|))$.

Note that if $(h_n)$ is a decreasing sequence of continuous functions on $\R_+$ converging to the characteristic function of $\{0\}$,  we have $\dim_\tau(\ker a)=\lim _n \tau(h_n(|a|))$.

Lück's method of handling Atiyah's problem can be understood in terms of traces via the following quite easy fact:

\begin{proposition}\label{propAt}
Let $A$ be a unital $C^*$-algebra, $a\in A$ and $s\in \R$. The map $\tau \mapsto \dim_\tau(\ker a)$ is continuous on $\Lambda_{a,s}$.
\begin{proof}
For $s'\in \R_+$, the set of traces $\Omega_{a,s'}=\{\tau;\ \tau(|a|)<s'\}$ is open. We only need to establish continuity on $\Lambda_{a,s}\cap \Omega_{a,s'}$. For $t\in \R_+$, put $\theta (t)=t-\ln _+(t)$. The function $\theta:\R_+\to \R$ is  continuous on $\R_+^*$, satisfies $\theta (t)>0$ for $t\ne 0$ and $\lim _{t\to 0^+} \theta(t)=+\infty$. Moreover, for every $\tau \in \Lambda_{a,s}\cap \Omega_{a,s'}$, we have $\tau(\theta (|a|))\le s'-s$.

The proposition is an immediate consequence of the following Lemma.
\end{proof}
\end{proposition}

\begin{lemme}
Let $\theta:\R_+\to \R_+$ be continuous on $\R_+^*$, and satisfy $\theta (t)>0$ for $t\ne 0$ and $\lim _{t\to 0^+} \theta(t)=+\infty$. Let $h_n:\R_+\to [0,1]$ be a decreasing sequence of continuous functions converging (pointwise) to $\chi_0$. Let $m\in \R_+$ and denote by $\Lambda_{a,m,\theta}$ the set of traces $\tau $ on $A$ such that $\tau(\theta (|a|))\le m$. Then the sequence $\Big(\tau(h_n(|a|))\Big)$ converges to $\dim_\tau(\ker a)$ uniformly on $\Lambda_{a,m,\theta}$.
\begin{proof}
The functions $v_n$ defined on $\R_+$ by $v_n(t)=\begin{cases} 
\frac{h_n(t)}{\theta(t)}& \hbox{if }\  t\ne 0\\
0& \hbox{if }\ t=0
\end{cases}$ are continuous. The sequence $(v_n)$ decreases to $0$; by Dini's theorem it converges uniformly to $0$ on $\Sp\,|a|$. For $\tau \in  \Lambda_{a,m,\theta}$, we have $0\le \tau(h_n(|a|))-\dim _\tau(\ker a)=\tau(v_n\theta(|a|))\le m\|v_n\|_\infty.$
\end{proof}
\end{lemme}

It follows that, if the group $\Gamma$ is either a direct limit, or a subgroup of an inverse limit (\emph{e.g.} a residually finite group), or a sofic group, then for any $a\in M_k(\Z [\Gamma])$ the von Neumann dimension of $\ker a$ can be computed as a limit of simpler terms.

\subsection{Algebraic coefficients}

We now see that proposition \ref{propAt} applies also for $a\in M_n(\C\Gamma)$ with algebraic coefficients (see \cite{DLMSY}). 

If $A$ is a unital $C^*$-algebra, $a\in A$ and $\tau$ is a trace on $A$, we define the rank $\mathrm{rk}_{\tau}(a)$ of $a$ to be the trace $\tau$ as the von Neumann dimension of the closure of the image of $a$; we of course have  $\mathrm{rk}_{\tau}(a)={\rm codim}_\tau(\ker a)=\tau(1)-\dim_\tau(\ker a)$. We will use the following Lemma.

\begin{lemme}\label{major}
 Let $A$ be a unital $C^*$-algebra, $a,b,c\in A$ with $a$ and $c$ invertible.  Then for any trace $\tau$ on $A$ we have $\tau(\ln_+ |abc|)\le \ln(\|a\|\|c\|)\mathrm{rk}_{\tau}(b)+\tau (\ln_+|b|)$.
\begin{proof}
We first prove this inequality when $c=1$: we prove, for invertible $a$,
 $$\tau(\ln_+ |ab|)\le \ln(\|a\|)\mathrm{rk}_{\tau}(b)+\tau (\ln_+|b|).\eqno(**)$$

We may replace $A$ by $\pi_\tau(A)''$ and thus assume $A$ is a von Neumann algebra and $\tau$ a faithful normal trace. Let $b=u|b|$ be the polar decomposition of $b$. Let $p=u^*u$ be the domain projection of $b$ and $A_p=pAp$. 
Let $\tau_p$ be the restriction of $\tau $ to $A_p$.

Now $|b|$ and $|au|$ are injective elements of $A_p$ so that the Fuglede-Kadison determinant $\Delta_{\tau_p}$ and Lück's modified determinant $\Delta_{\tau_p}^+$ coincide on them. We therefore have $\tau(\ln_+(|au||b|))=\tau(\ln_+(|au|))+\tau(\ln_+|b|))$. Finally, $\ln_+(|au|)\le \ln\|a\|p$, and $\tau (p)=\mathrm{rk}_{\tau}(b)$.

To end, since $\tau(\ln_+ |x|)=\tau(\ln_+|x^*|)$ we find (using (**))  $$\tau(\ln_+ |bc|)\le \ln(\|c\|)\mathrm{rk}_{\tau}(b^*)+\tau (\ln_+|b|).$$

Replacing $b$ by $bc$ in (**) and noting that $\mathrm{rk}_{\tau}(b)=\mathrm{rk}_{\tau}(bc)$, we find: $$\tau(\ln_+ |abc|)\le \ln(\|a\|)\mathrm{rk}_{\tau}(bc)+\tau (\ln_+|bc|)\le \ln(\|a\|\|c\|)\mathrm{rk}_{\tau}(b)+\tau (\ln_+|b|).$$
\end{proof}
\end{lemme}

\begin{proposition}\label{propalgebraic} Let  $a\in M_n(\overline\Q \F_\infty)$.
\begin{enumerate}\renewcommand \theenumi{\alph{enumi}}
\renewcommand \labelenumi{\rm {\theenumi})}
\item  There exists a constant $m$ such that for any tracial state $\tau\in \Lambda_{\F_\infty}$ (\ie a tracial state on $C^*(\F_\infty)$ satisfying Lück's property) we have $(\tau \otimes \Tr_n)(\ln_+|a|)\ge m$.\label{propalgebraica}
\item The map $\tau \mapsto \dim_\tau(\ker a)$ is continuous on $\Lambda_{\F_\infty}$.\label{propalgebraicb}
\end{enumerate}
\begin{proof}
Note that (\ref{propalgebraicb}) is an immediate consequence of (\ref{propalgebraica}) and prop. \ref{propAt}.

We prove (\ref{propalgebraica}). Let $K$ be a finite extension of $\Q$ containing the coefficients of $a$, \ie such that $a\in M_d(K \F_\infty)$. 
 
 Choosing a $\Q$ basis of $K$,  we obtain an embedding $i:K\to M_d(\Q)$ (where $d$ is the dimension of $K$ over $\Q$).
 
 Giving all the embeddings of $K$ to $\C$, we obtain an embedding $j=(j_1,\ldots,j_d):K\to \C^d$. We will assume that the given embedding $K\subset \C$ is $j_1$.
 
 These two embeddings are conjugate in $M_d(\C)$: There exists an invertible matrix $c\in M_d(\C)$ such that, for $x\in K$, the matrix $c\,i(x)c^{-1}$ is the diagonal matrix with coefficients $j_\ell (x)$.
 
 Write then $i(a)=k^{-1}b$ where $b\in M_{dn}(\Z \F_\infty)$ and $k\in \N^*$. 
 
 For every $\tau \in\Lambda_{\F_\infty}$, we have\begin{itemize}
\item $\tau (\ln_+(b))\ge 0$ (by definition of $\Lambda_{\F_\infty}$);
\item it follows that $\tau (\ln_+(i(a)))\ge -nd \ln k$;
\item using lemma \ref{major}, we find $\sum_{\ell=1}^d\tau (\ln_+(j_\ell (a)))\ge -nd \ln k-nd\ln(\|c\|\|c^{-1}\|)$.
\item On the other hand for all $\ell$, we have $\tau (\ln_+(j_\ell (a)))\le n\max(0,\ln \|j_\ell (a)\|)$;
\item we find $\tau (\ln_+(a))=\tau (\ln_+(j_1 (a)))\ge -nd \ln (k\|c\|\|c^{-1}\|)-n\sum_{\ell=2}^d\max(0,\ln \|j_\ell (a)\|)$.
\qedhere
\end{itemize}
\end{proof}
\end{proposition}

Generalizing a result of \cite{DLMSY}, we find:

\begin{proposition} Let $a\in M_n(\overline\Q \F_\infty)$. For any sofic trace of $C^*(\F_\infty)$, the dimension $\dim_\tau (\ker a)$ does not depend on the embedding $\overline\Q\subset \C$.

\begin{proof}
Denote by $\Sigma_{a}$ the set of tracial states on $C^*(\F_\infty)$ satisfying Lück's property and for which $\dim _\tau(\ker j(a))$ does not depend on the embedding $j:\overline\Q\to \C$. We wish to prove that every sofic trace is in $\Sigma_{a}$.

It follows from proposition \ref{propalgebraic} (applied to $a$ and $j(a)$ where $j$ is another inclusion of $\overline \Q$ in $\C$) that the set $\Sigma_{a}$ is closed in $\Lambda_{\F_\infty}$ and therefore in $\cT_{\F_\infty}$.

Let $q:\F_\infty\to \gS_k$ be a morphism and $\tau_q$ the corresponding trace on $C^*(\F_\infty)$. Denote still by $q$ the corresponding morphism $q:M_n(\overline \Q\F_\infty)\to M_{kn}(\overline \Q)$. For any embedding $j:\overline \Q\to \C$, we have $j\circ q=q\circ j:M_n(\overline \Q\F_\infty)\to M_{kn}(\C)$, so that we have $$\dim _{\tau_q}(\ker j( a))=\frac{1}{k}\dim_\C j\circ q( a)=\frac{1}{k}\dim_{\overline \Q}  q( a).$$
It is independent of $j$. 

In other words, $\Sigma_{a}$ is closed and contains all permutational traces. Therefore $\Sigma_{a}$ contains the closure of permutational traces: the sofic traces.
\end{proof}
\end{proposition}

We immediately find the following:

\begin{corollaire}
\label{corindep} Let $\Gamma$ be a sofic group and $a\in M_n(\overline\Q \Gamma)$, then the von Neumann dimension (with respect to the group trace of $\Gamma$) of $\ker a$ does not depend on the embedding $\overline\Q\subset \C$. \hfill $\square$
\end{corollaire}

A particular case of this corollary is proved by A. Thom \cite[theorem 4.3.(ii)]{Thom}. Note that Thom's proof actually establishes this result.

\bibliography{BiblioSG.bib} 
\bibliographystyle{amsplain}

\end{document}